\documentclass[12pt,leqno,twoside]{amsart}
\usepackage{amssymb}
\usepackage{latexsym}

\newtheorem{theorem}{Theorem}[section]

\newtheorem{proposition}[theorem]{Proposition}
\theoremstyle{definition}
\newtheorem{definition}[theorem]{Definition}
\theoremstyle{remark}
\newtheorem{remark}[theorem]{\sc Remark}
\newtheorem{example}[theorem]{\sc Example}

\newcommand{\reg}{{\rm{reg}}}
\renewcommand{\r}{{\rm{rHd\hspace{2pt}}}}

\newcommand{\Eu}{{\rm{Eu}}}

\newcommand{\rank}{{\rm{rank\hspace{2pt}}}}

\newcommand{\grad}{\mathop{\rm{grad}}\nolimits}

\renewcommand{\d}{{\rm{d}}}

\newcommand{\m}{\setminus}

\newcommand{\cW}{{\mathcal W}}

\newcommand{\bC}{{\mathbb C}}

\newcommand{\bv}{{\bf{v}}}

\begin{document}

\title[Milnor numbers and Euler obstruction]
{Milnor numbers and Euler obstruction}

\author{\sc Jos\'e Seade}

\address{J.S. and A.V.: Instituto de Matem\'aticas, Unidad Cuernavaca, Universidad Nacional 
Aut\'onoma de M\'exico, Apartado postal 273-3, C.P. 62210, Cuernavaca, Morelos, M\'exico.}

\email{jseade@matem.unam.mx}\email{alberto@matcuer.unam.mx}

\author{\sc Mihai Tib\u ar}

\address{M.T.:  Math\' ematiques, UMR 8524 CNRS,
Universit\'e de Lille 1, \  59655 Villeneuve d'Ascq, France.}

\email{tibar@math.univ-lille1.fr}

\author{\sc Alberto Verjovsky}

\thanks{Partially supported by CNRS-CONACYT (12409) Cooperation Program. 
The first and third named authors partially supported by
CONACYT grant G36357-E and DGPA (UNAM) grant IN 101 401}






\begin{abstract}
Using a geometric approach, we determine the relations between
the local Euler obstruction $\Eu_f$ of a holomorphic function $f$ and several
generalizations of the Milnor
number for functions on singular spaces.

\end{abstract}
%
\maketitle

\setcounter{section}{0}
\section{Introduction}\label{intro}
In the case of a nonsingular germ $(X,x_0)$ and a function $f$ with an
isolated critical point at $x_0$, the following three invariants
coincide (for (c), up to sign):
 
\begin{enumerate}
\item the Milnor number of $f$ at $x_0$, denoted $\mu(f)$;
\item  the number of Morse points in a Morsification of $f$;
\item  the Poincar\'e-Hopf index of $\overline{\grad f}$ at $x_0$;
\end{enumerate}
  
This fact is essentially due to Milnor's work in the late sixties
\cite{Mi}.
  There exist extensions of all these invariants to the case when
  $(X,x_0)$ is a singular germ, but they do not coincide in general.
One of the extensions of (c) is  the Euler obstruction of $f$ at $x_0$,
  denoted $\Eu_f(X,x_0)$. This was introduced in \cite{BMPS}; roughly, it is
  the obstruction to extending the conjugate of the gradient of the 
function $f$ as a section of the Nash bundle of $(X,x_0)$. 
It  measures how far the local Euler obstruction is
from satisfying the local Euler condition with respect 
to $f$ in bivariant theory.  
  It is then natural to compare
  $\Eu_f(X,x_0)$ to the Milnor number of $f$ in the case of a
  singular germ $(X,x_0)$. This has been also a question raised in \cite{BMPS}.

The main idea of this paper is that, for singular $X$, the Euler
 obstruction $\Eu_f(X,x_0)$ is most closely related to (b). We use the
 homological version of the
bouquet theorem for the Milnor fiber given in \cite{Ti}, which
relates the contributions in the bouquet to the
number of Morse points. Through this relation,
one may compare $\Eu_f(X,x_0)$ to the highest Betti number of the Milnor fiber
 of $f$. In case $X$ has {\em Milnor's property}, the comparison is optimal
 and yields a general inequality, see \S \ref{ss:milnor}.   
We further compare $\Eu_f(X,x_0)$ with two different generalizations
of the Milnor number for functions with isolated singularity on
singular spaces, one due to
\cite{Le-isol}, the other to
\cite{Go,MS} for curve singularities and to \cite {IS} for functions 
on isolated complete intersection germs in general.
 In case when the germ $(X,x_0)$ is an
isolated complete intersection singularity, we use in addition the
GSV-index of vector fields \cite{GSV} to completely determine the relations
between $\Eu_f(X,x_0)$, the Milnor number of $f$ and the GSV-index attached to 
$f$. 




\section{Euler obstruction and Morsification of
functions}\label{pre}

Let $(X,x_0)$ denote the germ at some point $x_0$ of a reduced pure
dimensional complex analytic space embedded in $\bC^N$, for some
$N$. Consider a Whitney stratification $\cW$ of some representative of
$X$. Let $W_0$ be the stratum containing $x_0$ and let $W_1, \ldots
W_q$ be the finitely many strata of $X$ having $x_0$ in their closure,
other than $W_0$. Let also $F: (\bC^N,x_0) \to (\bC, 0)$ denote some extension of $f$.

\begin{definition}\label{d:morse} (Lazzeri '73, Benedetti '77, Pignoni
'79, Goresky-MacPherson '83 \cite[p.52]{GM}.) One says that $f: (X,
x_0) \to \bC$ is a {\em general function} if $\d F_{x_0}$ does not
vanish on any limit of tangent spaces to $W_i$, $\forall i\not= 0$,
and to $W_0\setminus \{ x_0\}$.  One says that $f: (X, x_0) \to \bC$
is a {\em stratified Morse function germ} if: $\dim W_0 \ge 1$, $f$ is
general with respect to the strata $W_i$, $i\not= 0$ and the
restriction $f_{|_ {W_0}}: W_0 \to \bC$ has a Morse point at $x_0$.
\end{definition}

Let us recall some definitions and notations from \cite{BMPS}. The complex conjugate of the gradient of the extension $F$ projects to the
tangent spaces of the strata of $X$ into a vector field, which may not
be continuous. One can make it continuous by ``tempering'' it in the
neighborhoods of ``smaller'' strata. One gets a well-defined
continuous stratified vector field, up to stratified homotopy, which
we denote by $\grad_X f$. We shall call it briefly {\it the gradient
vector field}.

 If $f$ is a function on $X$ with an isolated singularity at $x_0$,
 with respect to the stratification $\cW$, then $\grad_X f$ has an
 isolated zero at $x_0$.  If $\nu: \widetilde X \to X$ is the Nash
 blow-up of $X$ and $\widetilde T$ is the Nash bundle over $\widetilde
 X$, then $\grad_X f$ lifts canonically to a never-zero section
 $\widetilde{\grad_X f}$ of $\widetilde T$ restricted to $\widetilde X
\cap \nu^{-1}(X \cap S_{\varepsilon})$, where $S_{\varepsilon}$ is a
small enough sphere around $x_0$, given by Milnor's result
\cite[Cor. 2.8]{Mi}.  Following \cite{BMPS}, the obstruction to extend
$\widetilde{\grad_X f}$ without zeros throughout $\nu^{-1}(X \cap B_{\varepsilon})$ is denoted by $\Eu_f(X,x_0)$ and is
called the {\em local Euler obstruction of $f$}.

\begin{example}\label{e:simple} 
If the germ $(X,x_0)$ is nonsingular,
then its Nash blow-up can be identified to $X$ itself, the Nash bundle is the usual
tangent bundle of $X$ and $\Eu_f(X,x_0)$ is, by definition, the
Poincar\'e-Hopf index of $\grad_X f$ at $x_0$.  From \cite[Th.7.2]{Mi} one deduces:
$\Eu_f(X,x_0)\,=\,(-1)^{\dim_\bC X} \mu$, where $\mu$ is the Milnor
number of $f$. It is also easy to prove (see \cite{BLS,BMPS}) that if
$(X,x_0)$ is any singular space but $f$ is a general function germ at
$x_0$, then the obstruction $\Eu_f(X,x_0)$ is zero.
\end{example}

 We claim that a natural way to study $\Eu_f$ is to split it according to a Morsification of $f$. We prove the following  general
formula for holomorphic germs with isolated singularity:
\begin{proposition}\label{p:euler} 
Let $f :(X,x_0)\to (\bC, 0)$ be a
holomorphic function with isolated singularity at $x_0$. Then
\[ \Eu_f(X,x_0) = (-1)^{\dim_\bC X} \alpha_q,\] where $\alpha_q$ is
the number of Morse points on $W_q = X_\reg$ in a generic deformation
of $f$.
\end{proposition}
\begin{proof} We Morsify the function $f$, i.e. we consider a small
analytic deformation $f_\lambda$ of $f$ such that $f_\lambda$ only has
stratified Morse points within the ball $B$ and it is general in a
small neighborhood of $x_0$. (See, for instance, the Morsification
Theorem 2.2 in \cite{Le}.)

Since $f_\lambda$ is a deformation of $f$, it follows that $\grad_X f$
is homotopic to $\grad_X f_\lambda$ over the sphere $X \cap \partial
B$, so the obstructions to extend their lifts to $\nu^{-1}(X \cap B)$
without zeros are equal.

On the other hand, the obstruction corresponding to $\grad_X
f_\lambda$ is also equal to the sum of local obstructions due to the
Morse points of $f_\lambda$.  Lemma 4.1 of \cite{STV} shows that the local
obstruction at a stratified Morse point is zero if the point lies in a
lower dimensional stratum.  So the points that only count are the
Morse points on the stratum $X_\reg$ and, at such a point, the
obstruction is $(-1)^{\dim_\bC X}$, as explain above in Example \ref{e:simple}.
\end{proof}

\begin{remark}\label{n:morse} 
The Euler obstruction is defined via the Nash blow-up and the latter only
takes into account the closure of the tangent bundle over the regular
part $X_\reg$.  Since the other strata are not counting in the Nash
blow-up, it is natural that they do not count for $\Eu_f(X,x_0)$
neither. The number $\alpha_q$ does not depend on
the chosen Morsification, by a trivial connectedness argument.
We refer to \cite{STV} for more about $\alpha_q$ and other invariants
of this type, which enter in a formula for the {\em global Euler obstruction}
of an affine variety $Y\subset \bC^N$.
\end{remark}
\begin{remark}\label{n:morse2}
The number $\alpha_q$ may be interpreted as the intersection number within
$T^*\bC^N$ between $\d F$ and the conormal $T^*_{X_\reg}$. Therefore
 our Proposition \ref{p:euler}
 may be compared to \cite[Corollary 5.4]{BMPS}, which is proved by using
different methods.
 J. Sch\"urmann informed us that such a result can also be obtained
using the techniques of \cite{Sch-book}.
\end{remark}


\section{Milnor numbers}\label{remarks}
\subsection{L\^e's Milnor number}\label{ss:milnor}
 L\^e D.T. \cite{Le-isol} proved that for a function $f$ with an isolated
singularity at $x_0 \in X$ (in the stratified sense) one has a Milnor fibration. He
pointed out that, under certain conditions, the space $X$ has
``Milnor's property'' in homology (which means that the reduced homology of
the Milnor fiber of $f$ is concentrated in dimension $\dim X
-1$). Then {\it the Milnor number} $\mu(f)$ is well defined as the
rank of this homology group.  By L\^e's results \cite{Le-isol}, Milnor's property is satisfied for instance if $(X,x_0)$ is a complete intersection
(not necessarily isolated!) or, more generally, if $\r (X,x_0)\ge \dim (X,x_0)$,
where $\r (X,x_0)$ denotes the {\em rectified homology depth} of $(X,x_0)$,
see \cite{Le-isol} for its definition originating in Grothendieck's work.

 To compare $\mu(f)$ with $\Eu_f(X,x_0)$ we use the general bouquet
 theorem for the Milnor fiber in its homological version. Let $M_f$
 and $M_l$ denote the Milnor fiber of $f$ and of a general function
 $l$. Let 
$f : (X,x_0) \to (\bC,0)$ be a function with stratified isolated
singularity and let $\Lambda$ be the set of stratified Morse points in
some chosen Morsification of $f$ (by convention $x_0\not\in
\Lambda$). Then by \cite[pp.228-229 and Bouquet Theorem]{Ti} we have:
 \begin{equation}\label{eq:bouquet}
 \tilde H_* (M_f) \simeq \tilde H_* (M_l) \oplus
\oplus_{i\in \Lambda}  H_{* - k_i +1} ( C (F_i), F_i)\,
\end{equation}
where, for $a_i\in \Lambda$, $F_i$ denotes the complex link of the
stratum to which $a_i$ belongs, $k_i$ is the dimension of this stratum
and $C (F_i)$ denotes the cone over $F_i$.

In particular, if the germ $(X,x_0)$ is a complete intersection (more
generally, if $\r (X,x_0)$ $\ge \dim (X,x_0)$), then: $\mu(f) = \mu(l)
+ \sum_{a_i\in \Lambda} \mu_i$, where $\mu_i := \rank H_{\dim X - k_i}
(C (F_i), F_i)$. 
This result shows that the Milnor number $\mu(f)$ gathers information
from all stratified Morse points, whereas $\Eu_f(X,x_0)$ is, up to
sign, the number $ \alpha_q = \#\Lambda_0$, where $\Lambda_0$ denotes
the set of Morse points occurring on $X_\reg$ (see Proposition \ref{p:euler}
above). Notice that we have $\Lambda_0 \subset \Lambda$, $\mu(l) \ge 0$, $\mu_i = 1$ if $i\in
\Lambda_0$ and $\mu_i \ge 0$ if $i\in \Lambda \setminus \Lambda_0$.  We
therefore get the general inequality, whenever the space $X$ has Milnor's property (e.g. when $(X,x_0)$ is a complete intersection,
not necessarily with isolated singularities), and therefore the
Milnor-L\^e number is well defined:
\begin{equation}\label{eq:mu}
\mu(f) \ge (-1)^{\dim X} \Eu_f(X,x_0).
\end{equation}
 In case $(X,x_0)$ is an isolated complete intersection singularity
 (ICIS for short), from the above discussion on (\ref{eq:bouquet}) we
 get the equality:
\begin{equation}\label{eq:mu2} 
\Eu_f (X,x_0) = (-1)^{\dim X} [\mu(f) - \mu(l)]. 
\end{equation}
In the ICIS case, (\ref{eq:mu2}) also shows that the inequality
(\ref{eq:mu}) is strict whenever $X$ is actually singular. This is so 
since $\mu(l) > 0$, which can be proved inductively using 
Looijenga's results \cite{Lo}.

\subsection{Another Milnor number}
A different generalization of the Milnor number is due to V. Goryunov
\cite{Go}, D. Mond and D. van Straten \cite{MS}. This is originally
 defined for functions on curve singularities $X \subset \bC^N$, and
 we refer to \cite[p.178]{MS} for the precise definition.  This
 number is preserved under simultaneous deformations of both the space
 $X$ and the function $f$. Thus, if the curve singularity $(X,x_0)$ is
 an ICIS, defined by some application $g\colon (\bC^N, x_0) \to
 (\bC^p, 0)$ on an open set in $\bC^N$, and $F$ is an extension of $f$
 to the ambient space, then $\mu_G(f)$ counts the number of critical
 points (with their multiplicities) of the restriction of $F$ to a
 Milnor fiber of $g$, say $X_t = g^{-1}(t)$ for some regular value $t$
 of $g$. This is equivalent to saying that $\mu_G(f)$ is the
 Poincar\'e-Hopf index of the gradient of the restriction
 $F_{|_{X_t}}$. In other words, this is saying that $\mu_G(f)$ is the
 {\em GSV-index} of the gradient vector field of $f$ on $X$. We recall
 that the GSV-index of a vector field $\bv$ on $(X,x_0)$, defined in
 \cite{GSV,SS}, equals the
 Poincar\'e-Hopf index of an extension of $\bv$ to the Milnor fiber
 $X_t$.
 
As noted in the introduction to \cite{BMPS}, this definition of
$\mu_G(f)$ makes sense in all dimensions and one may generalize
$\mu_G$ as follows.  Given an ICIS $(X,x_0)$ and a function $f$ on it
with an isolated singularity at $x_0$, we denote by $\nabla_X f$ the
gradient vector field of $f$ (not the conjugate of the gradient as we
did for defining $\Eu_f(X,x_0)$). 
Thus we may define $\mu_G(f)$ as the
GSV-index of $\nabla_X f$ at $x_0$.
We notice that this invariant is precisely  the 
{\em virtual multiplicity} at $x_0$ of the function $f$ on $X$
 introduced by Izawa and Suwa in
\cite{IS} and denoted   $\tilde m(f;x_0)$. This multiplicity is by 
definition
the localization at $x_0$ of the top Chern class of the 
virtual cotangent bundle $T^*(X)$ of $X$  defined by 
the differential of $f$, which is non-zero on $X \m \{x_0\}$ by 
hypothesis. This invariant has the advantage of being defined even
if the singular set of $X$ is non-isolated and it is related to global
properties of the variety 
(we refer to \cite{IS} for details). 
This coincides with the index
of the 1-form $dg$ defined in \cite{EG} and it is similar to the interpretation of the GSV index 
of  vector fields given in \cite{LSS} as a  
localization  of the top Chern class of the virtual tangent bundle.

One can easily find the relation between $\mu_G(f)$ and
$\mu(f)$ in case $X$ is an ICIS. The proof can be found for instance
in \cite{Lo}. Let $\mu(X,x_0)$ be the Milnor number of  the
 ICIS $(X, x_0)$ and let $f$
be some function with isolated singularity on $(X, x_0)$. Then:
\[ \mu_G(f) = \mu(f) + \mu(X,x_0).\]

Using (\ref{eq:mu2}) we get:
\begin{equation}\label{eq:mu3} 
\Eu_f (X,x_0) = (-1)^{\dim X} [\mu_G(f) - \mu_G(l)].
\end{equation}
These equalities completely determine the relation
between $\Eu_f(X,x_0)$, the GSV-index and the Milnor number of $f$, in
terms of the Milnor number of the ICIS $(X,x_0)$.

\section{Further remarks}
 
It is proved in \cite{BMPS}, using
\cite{BLS}, that one has:
\begin{equation}\label{eq:releuler} \Eu_f(X,x_0) =
\sum_{i=0}^q[\chi(M(l,x_0) \cap W_i) - \chi(M(f,x_0) \cap W_i)] \cdot
\Eu_X(W_i),
\end{equation} 
where $M(f,x_0)$ and $M(l,x_0)$ denote representatives
of the Milnor fibers of $f$ and of the generic linear function $l$, respectively.  Combining this
relation with Proposition \ref{p:euler} one gets:
\[ \sum_{i=0}^q[\chi(M(l,x_0) \cap W_i)
- \chi(M(f,x_0) \cap W_i)] \cdot
\Eu_X(W_i) = (-1)^{\dim_\bC X} \alpha_q.
\]

\begin{example}\label{e:morse} Let $X = \{ x^2 - y^2 = 0\} \times \bC \subset \bC^3$ and $f$ be the restriction to $X$ of the function 
$(x,y,z) \mapsto x+ 2y + z^2$.  Take $x_0 := (0,0,0)$ and take as
general linear function $l$ the restriction to $X$ of the projection
$(x,y,z) \mapsto z$.  Then $X$ has two strata: $W_0 = $ the $z$-axis,
$W_1 = X\m \{ x=y=0\}$.  We compute $\Eu_f(X,x_0)$ from the relation
(\ref{eq:releuler}).

 First, $M(l,x_0)\cap W_0$ is one point and $M(f,x_0)\cap W_0$ is two
 points. Next, $M(l,x_0)\cap W_1$ is the disjoint union of two copies
 of $\bC^*$ and $M(f,x_0)\cap W_1$ is the disjoint union of two copies
 of $\bC^{**}$, where $\bC^*$ is $\bC$ minus a point and $\bC^{**}$ is
 $\bC$ minus two points.  Then formula (\ref{eq:releuler}) gives: $
 \Eu_f(X,x_0)= (1-2)\cdot \Eu(X,x_0) + (0-(-2)) \cdot 1$.
 
We have $\Eu(X,x_0)=\Eu(X\cap \{l=0\}, x_0)$. Next $\Eu(X\cap \{l=0\},
x_0)$ is just the Euler characteristic of the complex link of the
slice $X\cap \{l=0\} = \{ x^2- y^2=0\}$. This complex link is two
points, so $\Eu(X\cap \{l=0\}, x_0) = 2$.  We therefore get
$\Eu_f(X,x_0)= 0$.
\end{example}



\end{document}